\documentclass{amsart}

\usepackage{breqn}
\usepackage{mathtools}
\usepackage{accents}

\begin{document}

\title[Arithmetic Euler top]{Arithmetic Euler top}
\bigskip

\def \h{\hat{\ }}
\def \cO{\mathcal O}
\def \ra{\rightarrow}
\def \bZ{{\mathbb Z}}
\def \cP{\mathcal V}
\def \cH{{\mathcal H}}
\def \cB{{\mathcal B}}
\def \d{\delta}
\def \cC{{\mathcal C}}
\def \jor{\text{jor}}

\newtheorem{THM}{{\!}}[section]
\newtheorem{THMX}{{\!}}
\renewcommand{\theTHMX}{}
\newtheorem{theorem}{Theorem}[section]
\newtheorem{corollary}[theorem]{Corollary}
\newtheorem{lemma}[theorem]{Lemma}
\newtheorem{proposition}[theorem]{Proposition}
\newtheorem{thm}[theorem]{Theorem}
\theoremstyle{definition}
\newtheorem{definition}[theorem]{Definition}
\theoremstyle{remark}
\newtheorem{remark}[theorem]{Remark}
\newtheorem{example}[theorem]{\bf Example}
\numberwithin{equation}{section}
\address{Department of Mathematics and Statistics\\University of New Mexico \\
  Albuquerque, NM 87131, USA\\ 
   Department of Mathematics and
  Statistics\\Boston University \\Boston, MA 02215, USA}
\email{buium@math.unm.edu \\ep@bu.edu } 
\maketitle

\bigskip

\medskip
\centerline{\bf Alexandru Buium and Emma Previato}
\bigskip

\begin{abstract}
The theory of differential equations has an arithmetic analogue in which derivatives of functions are replaced by Fermat quotients of numbers. Many classical differential equations (Riccati, Weierstrass, Painlev\'{e}, etc.) were previously shown to possess arithmetic analogues. The paper introduces an arithmetic analogue of the Euler differential equations for the rigid body.
\end{abstract}

\section{Introduction}

The theory of differential equations has an arithmetic analogue in which derivatives  of functions are replaced by Fermat quotients of numbers. This arithmetic analogue  was introduced in \cite{char}; for an exposition of part of the resulting theory  we refer to \cite{book}.
The present paper morally fits into the theory developped in \cite{char,book}. However our paper is written so as to be entirely self-contained and, in particular, it is independent of \cite{char, book}; the few facts we need from \cite{char,book}  will be quickly reviewed here.

Many remarkable classical differential equations have arithmetic analogues. Examples of such classical differential equations are: the Riccati equation \cite{adel3}, the Weierstrass equation \cite{char},  the Painlev\'{e} VI equation \cite{char, BYM}, 
Schwarzian equations (satisfied by modular forms) \cite{char}, 
and linear differential equations corresponding to special connections in Riemannian geometry (such as Chern, Levi-Civita, etc.) \cite{adel2, adel3, curvature1, curvature2}.

The purpose of the present paper is to develop an arithmetic analogue of the Euler differential equations for the rigid body (the Euler top). This is a system of $3$ ordinary (non-linear) differential equations in $3$ variables which is one of the simplest examples of   {\it algebraically completely integrable} systems \cite{moerbeke}. As such its flow on $3$-space, referred to in what follows as the {\it classical Euler flow}, can be viewed as a derivation $\d$ on the polynomial ring ${\mathbb C}[x_1,x_2,x_3]$ given by 
\begin{equation}
\label{classical EF}
\d=(a_2-a_3)x_2x_3\frac{\partial}{\partial x_1}+(a_3-a_1)x_3x_1\frac{\partial}{\partial x_2}+(a_1-a_2)x_1x_2\frac{\partial}{\partial x_3},\end{equation}
where $a_1,a_2,a_3\in {\mathbb C}$ are distinct complex numbers. This flow is trivially seen to have
 $2$ independent  prime integrals 
\begin{equation}
\label{HH}
H_1=\sum_{i=1}^3a_ix_i^2,\ \ \ H_2=\sum_{i=1}^3x_i^2,\end{equation}
in the sense that
\begin{equation}
\label{HH0}
 \d H_1=\d H_2=0.\end{equation}
For generic $c=(c_1,c_2)\in {\mathbb C}^2$, the loci $E_c$ in $3$-space, given by 
\begin{equation}
\label{loci}
H_1=c_1,\ \ \ H_2=c_2,\end{equation}
  are affine elliptic curves. We refer to $E_c$ as the {\it level sets} of $H_1,H_2$.  Then  the  classical Euler flow $\d$ is  ``linearized" when restricted to these level sets $E_c$ in the sense that, if one denotes by $\d_c$ the action of $\d$ as Lie derivative on the $1$-forms on $E_c$ and if $\omega_c$ is the ``canonical" invariant $1$-form on $E_c$, then
  \begin{equation}
  \label{Lie}
  \d_c\omega_c=0.
  \end{equation}
  
  \medskip
 
   Here is  our arithmetic analogue of the above. 
      
   Let $A$ be a complete discrete valuation ring with maximal ideal generated by an odd prime $p$ and residue field ${\mathbb F}=A/pA$ algebraic over the prime field. Let $a_1,a_2,a_3\in A$ be distinct mod $p$ and consider again $H_1,H_2$ as in \ref{HH}. 
   Recall from \cite{char} that a {\it $p$-derivation} on a ring $B$ (in which $p$ is a non-zero divisor) is a map $\d:B\ra B$ such that the map $\phi:B\ra B$, $\phi(b):=b^p+p\d b$, is a ring homomorphism; $\phi$ is then referred to as the {\it Frobenius lift} attached to $\d$. 
   Of course, $A$ has a unique $p$-derivation; any $p$-derivation on an $A$-algebra will be tacitly  assumed  compatible with the one on $A$.
   We also denote by $\widehat{\ }$ the $p$-adic completion of rings or schemes.
   We will prove:
   
   \begin{thm}
   \label{main1}
   There exists an open set $X$ of the affine $3$-space over $A$, a $p$-derivation $\d$ on  $\cO(\widehat{X})$, with attached Frobenius lift $\phi$, and  an invertible function $\lambda\in \cO(X)^{\times}$ such that:
   
   1) The following equalities hold in $\cO(\widehat{X})$:
   $$\d H_1=\d H_2=0;$$
   
   2) For any $c\in A^2$  in a suitable Zariski open set and satisfying $\d c=0$ we have
   $$\frac{\phi_c^*}{p}\omega_c\equiv \lambda(c)\cdot \omega_c\ \ \ \text{mod}\ \ \ p,$$
   on $\widehat{E}_c\cap \widehat{X}$, where $E_c$ are the level sets of $H_1,H_2$, the map $\phi_c$ is the restriction of $\phi$ to  $\widehat{E}_c\cap \widehat{X}$, and $\omega_c$ is the invariant $1$-form on $E_c$.
   \end{thm}
   
    Conditions 1 and 2 in our theorem are of course to be viewed as analogues of conditions \ref{HH0} and \ref{Lie} respectively. For a more precise statement of the above theorem we refer to Theorem \ref{linearization theorem} in the text; this more precise statement will explicitly give $X$ and $\lambda$.
   As we shall see the function $\lambda$ will be the ``inverse of the Hasse invariant" for a family of plane quartics isogenous to the family of level sets of $H_1,H_2$.
   
   A $p$-derivation  $\d$  as  in  Theorem \ref{main1}     will be referred to as an {\it arithmetic Euler flow} on $\widehat{X}$. 
   Arithmetic Euler flows are, of course,  not unique. Remarkably, it will turn out that the set of arithmetic Euler flows mod $p$ is parameterized by the ``prime integrals mod $p$" for the {\it classical Euler flow} \ref{classical EF}. Indeed  formula \ref{classical EF} defines a derivation on the ${\mathbb F}$-vector space  $\cO(\overline{X})$, where $\overline{X}:=X \otimes {\mathbb F}$, and we have:
   
   \begin{thm}
   \label{main2}
   Assume ${\mathbb F}$ is infinite. Then
   the set of equivalences classes mod $p$ of arithmetic Euler flows on $\widehat{X}$ is a principal homogeneous space for the ${\mathbb F}$-vector space of all elements in $\cO(\overline{X})$ killed by  the derivation \ref{classical EF}.
   \end{thm}

    Cf. Theorem \ref{principal} in the text for a more precise statement in which $X$ and $\lambda$ are, again,  specified explicitly.
   
   The principal homogeneous space in Theorem \ref{main2} is non-empty by Theorem \ref{main1}.
   
It is reasonable to expect that, more generally, higher dimensional algebraically completely integrable systems might have similar arithmetic analogues; no straightforward generalization of our arguments here seems possible, however, at this point. 

The paper is organized as follows. Section 2 sets the geometric stage for the Euler equations (for both the classical and the arithmetic case). Section 3 reviews the {\it classical Euler flow}. Section 4 quickly reviews Fermat quotients ($p$-derivations) following \cite{char,book}. Section 5 reviews the Hasse invariant. In section 6 we state  our main results which are Theorems 
\ref{linearization theorem} and \ref{principal}. Section 7 will provide the proofs of our theorems.

\bigskip

{\bf Acknowledgments}. The authors are indebted to the IHES in Bures sur Yvette, where their collaboration on this projected started.
The first author was partially supported by the Simons Foundation (award 311773).
The second author expresses her sincerest
 gratitude to the Boston University
College of Arts and Sciences Associate Dean Stan Sclaroff and Assistant Dean
Richard Wright, for the travel support that made this research possible.

\section{Basic geometric setting}\label{section 2}

In this section we introduce the geometric objects that will be relevant throughout the paper.

Let $A$ be 
either a field or a discrete valuation ring and assume  $2$ is invertible in $A$.
We denote by $A^{\times}$ the group of invertible elements.
Let $x_1,x_2,x_3$ and $z_1,z_2$ be variables. Let $a_1,a_2,a_3\in A$ be such that 
\begin{equation}
\label{condition 1}
(a_1-a_2)(a_2-a_3)(a_3-a_1)\in A^{\times}.
\end{equation} 
We consider the quadratic forms $H_1,H_2\in A[x_1,x_2,x_3]$, 
\begin{equation}
\label{condition 2}
H_1:=\sum_{i=1}^3 a_ix_i^2,\ \ \ H_2:=\sum_{i=1}^3x_i^2.
\end{equation}
Also we consider the affine spaces 
\begin{equation}
\label{condition 3}
{\mathbb A}^2=Spec\ A[z_1,z_2],\ \ \ {\mathbb A}^3=Spec\ A[x_1,x_2,x_3]
\end{equation}
and the morphism
\begin{equation}
\label{condition 4}
{\mathcal H}:{\mathbb A}^3\ra {\mathbb A}^2
\end{equation}
defined by
\begin{equation}
z_1\mapsto H_1,\ \ \ z_2\mapsto H_2.\end{equation}
 For any $A$-point $c=(c_1,c_2)\in A^2={\mathbb A}^2(A)$ we denote by 
\begin{equation}
\label{definition of Ec}
E_c=Spec\ A[x_1,x_2,x_3]/(H_1-c_1,H_2-c_2)
\end{equation}
the fiber of ${\mathcal H}$ at $c$ and by $i_c:E_c\ra {\mathbb A}^3$ the inclusion. 

Consider the polynomial
\begin{equation}
 \label{definition of N}
 N(z_1,z_2)= \prod_{i=1}^3(z_1-a_iz_2)\in A[z_1,z_2].\end{equation}
 Assume, in addition, that 
  \begin{equation}
 \label{nondeg}
 N(c_1,c_2)\in A^{\times},\end{equation}
 equivalently that
 all $2\times 2$ minors of the matrix 
 \begin{equation}
 \label{the two by four matrix}
 \left(\begin{array}{cccc}
 a_ 1& a_2 & a_3 & c_1\\
 1 & 1 & 1 & c_2
 \end{array}\right)\end{equation}
 are invertible. Then, for $E_c^{ij}$ the locus in $E_c$ where $x_ix_j$ is invertible, we have that
 $$E_c=E^{12}_c\cup E^{23}_c\cup E_c^{31}$$
 and $E_c$ is smooth over $A$. 
 Moreover $E_c$ comes equipped with a global $1$-form
 given by
  \begin{equation}
  \label{the form on intersections of two quadrics}
  \omega_c=i_c^*\frac{dx_1}{(a_2-a_3)x_2x_3}=i_c^*\frac{dx_2}{(a_3-a_1)x_3x_1}=i_c^*\frac{dx_3}{(a_1-a_2)x_1x_2},\end{equation}
  which we will refer to as the {\it canonical} $1$-form on $E_c$.
If one considers  the  projective closure $\overline{E}_c$ of $E_c$ in the projective space $${\mathbb P}^3=Proj\ A[t_0,t_1,t_2,t_3],\ \ \ x_i=t_i/t_0,$$
then, by the arguments above applied to the other affine charts of
${\mathbb P}^3$ one gets that $\overline{E}_c$ is smooth and comes equipped with a $1$-form whose restriction to $E_c$ is $\omega_c$.
    
   Here are more facts, needed later,  about the geometry of the situation. Consider  two more indeterminates $x,y$,  and consider the  polynomial 
   $$F=F(z_1,z_2,x)\in A[z_1,z_2][x]$$
   defined by 
      \begin{equation}
     \label{quartic in one variable}
     F:=((a_2-a_3)x^2+z_1-a_2z_2)((a_3-a_1)x^2-z_1+a_1z_2).
     \end{equation}
     This is a weighted homogeneous polynomial of degree $4$ (a quartic) in $z_1,z_2,x$ with respect to the weights $(2,2,1)$.
   For any $c=(c_1,c_2)\in A^2$ set
   \begin{equation}
   E'_c:=Spec\ A[x,y]/(y^2-F(c_1,c_2,x)).
   \end{equation}
   Then we have a morphism
     \begin{equation}
     \label{isogeny}
    \pi: E_c\ra E'_c,\end{equation}
     given by
     $$ x\mapsto x_3,\ \ \ y\mapsto (a_1-a_2)x_1x_2.$$
     Indeed if we consider the ideal
     \begin{equation}
     \label{ideal I}
     I=(H_1-z_1, H_2-z_2)\subset A[x_1,x_2,x_3][z_1,z_2]
     \end{equation}
     we have the following congruences:
     \begin{equation}
     \label{congruences mod I}
     \begin{array}{rcll}
     (a_1-a_2)x_1^2 & \equiv & (a_2-a_3) x_3^2+ z_1-a_2z_2 & \text{mod}\ \ I\\
     \  & \  & \  & \  \\
      (a_1-a_2)x_2^2 & \equiv & (a_3-a_1) x_3^2- z_1+a_1z_2 & \text{mod}\ \ I
     \end{array}\end{equation}
     hence
     \begin{equation}
     \label{one more congruencies}
     (a_1-a_2)^2x_1^2x_2^2 \equiv F(x_3,z_1,z_2)\ \ \ \text{mod}\ \ I.\end{equation}
     Note that, under the assumption that $N(c_1,c_2)\in A^{\times}$, the discriminant  of $F$ is in $A^{\times}$ so $E'_c$ is smooth over $A$. Moreover,
      \begin{equation}
      \pi^*(\frac{dx}{y})=i_c^*\frac{dx_3}{(a_1-a_2)x_1x_2}=\omega_c.\end{equation}
     For $A$ a perfect  field and   $c_1,c_2$ satisfying $N(c_1,c_2)\neq 0$  we have that $E'_c$ is a smooth plane curve whose smooth projective model $\overline{E}'_c$ is an elliptic curve and hence \ref{isogeny} is induced by an isogeny of elliptic curves, 
     \begin{equation}
     \label{completed isogeny}
     \overline{E}_c\ra \overline{E}'_c.\end{equation}
    This isogeny has degree $2$, in particular it is separable.

\section{Classical Euler flow}

In this section we review the classical Euler equations for the rigid body (the Euler top).

With notation as in section \ref{section 2} assume, in the present section only, that $A$ is an algebraically closed field of characteristic $\neq 2$ (e.g., the complex field, which is the classical setting for the theory).  
Consider  the $A$-derivation $\d$ on the polynomial ring $A[x_1,x_2,x_3]$  given by
\medskip
 \begin{equation}
 \label{Euler system}\begin{array}{rcl}
 \d x_1 & = & (a_2-a_3)x_2x_3,\\
  \d x_2 & = & (a_3-a_1)x_3x_1,\\
   \d x_3 & = & (a_1-a_2)x_1x_2.\end{array}\end{equation}
   \medskip
   We we refer to the derivation $\d$ as   the {\it classical Euler flow}.
   
   For any $c=(c_1,c_2)\in A^2$ with $N(c_1,c_2)\neq 0$
   denote by $\d_c$ the derivation on 
   $$\cO(E_c)=A[x_1,x_2,x_3]/(H_1-c_1,H_2-c_2)$$ induced by the derivation $\d$ on $A[x_1,x_2,x_3]$. 
   
   Here are a couple of trivial facts about the derivation \ref{Euler system}
    that we collect in the form of the Theorem below; in this Theorem $\langle\ ,\ \rangle$ stands for the pairing between derivations and $1$-forms. 
   
   \begin{thm}
   
   \label{classical theorem}
  
  \ 
  
  1) The following equalities hold in $A[x_1,x_2,x_3]$:
  $$
  \d H_1=\d H_2=0.$$

   2) For any $c=(c_1,c_2)\in A^2$ with $N(c_1,c_2)\neq 0$ 
   the  derivation $\d_c$ on $\cO(E_c)$ and the canonical $1$-form  $\omega_c$ on $E_c$
 are dual to each other in the sense that
  $$
  \langle\d_c,\omega_c\rangle=1.
  $$\end{thm}
  
  \medskip
  
    \begin{remark}\

    1) Condition 1 above is interpreted as saying that $H_1,$ and $H_2$ are {\it prime integrals} for $\d$. 
    
    2) Condition 2 implies  that the derivation $\d_c$ on $E_c$ extends to a derivation of the structure sheaf of $\overline{E}_c$; morally this says that $\d_c$ is {\it linearized} on $E_c$. 
    
   3)  Let $\d_c$ be {\it any} $A$-derivation on $E_c$
   and let us still denote by $\d_c$ the action of the Lie derivative along $\d_c$ on $1$-forms on $E_c$.  Then we have
   \begin{equation}
   \label{monkey}
   \langle \d_c,\omega_c\rangle\in A\ \ \ \ \Leftrightarrow \ \ \ \ \d_c \omega_c=0.\end{equation}
    This follows from Cartan's formula for the Lie derivative (because $2$-forms on curves vanish) or, more directly, one can argue as follows. Set $\omega_c=fdx_3$  and $\d_c =g\frac{\partial}{\partial x_3}$.  
    So $\langle \d_c,\omega_c\rangle=fg$. Now,
$$  \begin{array}{rcl}
  \d_c\omega_c & = &\d_c\left(fdx_3\right)\\
  \  & \  & \  \\
  \  & = & (\d_c f) dx_3+f d(\d_c x_3)\\
   \  & \  & \  \\
  \  & = & 
g\cdot \frac{\partial f}{\partial x_3} \cdot dx_3+f\cdot dg\\
 \  & \  & \  \\
\  & = & gdf+fdg\\
 \  & \  & \  \\
\  & = & d(fg),\end{array}$$
which ends the proof of \ref{monkey}. In particular condition 2 in Theorem \ref{classical theorem} implies (and ``up to a scalar" is implied by)  $\d_c\omega_c=0$. The condition $\langle \d_c,\omega_c\rangle =1$ does not have a direct analogue in the arithmetic case; but, as we shall see,  the condition $\d_c\omega_c=0$ will.

      4) Assertions 1 and 2 in Theorem \ref{classical theorem} make the classical Euler flow  fit into the paradigm of  algebraic complete integrability \cite{moerbeke}. 
    \end{remark}
     
    \section{Review of the arithmetic setting}

Following \cite{char, book} we review, in this section, some of the basic definitions needed in the arithmetic setting.

 Let $p$ be an odd prime. The sign 
 $\widehat{ \ }$ will always mean $p$-adic completion of rings or schemes. So for a ring $B$ we have
 $$\widehat{B}=\lim_{\leftarrow} B/p^nB$$
while,  for an affine scheme $X=Spec\ B$, $\widehat{X}$ is the formal scheme \cite{hartshorne},
 $$\widehat{X}=Spf\ \widehat{B}.$$
 Recall that a Frobenius lift on a ring $B$ is a ring homomorphism $\phi=\phi^B:B\ra B$ whose reduction mod $p$ is the $p$-power Frobenius on $B/pB$. More generally, for a $B$-algebra $C$, a Frobenius lift from $B$ to $C$ is a ring homomorphism $B\ra C$ whose reduction mod $p$ is the structure map of the algebra composed with the $p$-power Frobenius. Given a ring $B$, in which $p$ is a non-zero divisor, a {\it $p$-derivation} on $B$ is a map of sets $\d=\d^B:B\ra B$ such that the map $\phi=\phi^B:B\ra B$ defined by 
 \begin{equation}
 \label{def of phi}
 \phi(b)=b^p+p\d b\end{equation}
  is a ring homomorphism; such a $\phi^B$ is, of course, a Frobenius lift and we say that $\d^B$ and $\phi^B$ are {\it attached} to each other. 
 More generally we may consider $p$-derivations from any ring $B$ into a $B$-algebra $C$ in which $p$ is a non-zero divisor; they are maps of sets such that \ref{def of phi} is a ring homomorphism. We shall still say that $\d$ and $\phi$ are attached to each other. If $X$ is a scheme and $\d$ is a $p$-derivation on $\cO(\widehat{X})$ (with $p$ a non-zero divisor in the latter) we also say that $\d=\d^X$ is a $p$-derivation on $\widehat{X}$; morally $\d^X$ should be viewed as an arithmetic analogue of a flow on $X$. The same convention will be in place for Frobenius lifts on formal schemes. The requirement above about $p$ being a nonzero divisor can be relaxed; cf. \cite{char, book}. This relaxation is necessary for the general theory but not for what we will need in the present paper.
 
  \medskip
  
  {\it Form now on, in this paper, $A$ will denote a complete discrete valuation ring  with maximal ideal generated by an odd prime $p$ such that the  residue field ${\mathbb F}:=A/pA$ is algebraic over the prime field ${\mathbb F}_p$.}
  
  \medskip
  
  Such an $A$ has a unique Frobenius lift and hence a unique $p$-derivation. From now all $p$-derivations, respectively Frobenius lifts, on $A$-algebras will be assumed compatible with the $p$-derivation, respectively the Frobenius lift, on $A$.

\section{Hasse invariant}

In this section we review some well known facts and trivial observations about the Hasse invariant.

\begin{definition}
Let $B$ be any ring and let $F\in B[x]$ be a polynomial in one variable $x$.
Define the {\it Hasse invariant} of $F$ to be the coefficient $A_{p-1}(F)\in B$ of $x^{p-1}$ in the polynomial
$F^{\frac{p-1}{2}}.$
\end{definition}

Let $C=Spec\ {\mathbb F}_p[x,y]/(y^2-F(x))$. 
We have the following well known facts:

\begin{lemma}
\label{well known}
Assume $B=\bZ_p$ is the ring of $p$-adic integers and $F\in \bZ_p[x]$ and set
$N_p=|C({\mathbb F}_p)|$.
 Then the following hold:
 
 1) If $\deg(F)=3$ then
$$ A_{p-1}(F)\equiv - N_p\ \ \ \text{mod}\ \ p.$$

2) If $\deg(F)=4$ and $F(x)=ax^4+\text{(lower terms)}$ then
$$ A_{p-1}(F)\equiv - N_p -  \left(\frac{a}{p}\right)\ \ \ \text{mod}\ \ p,$$
where $\left(\frac{a}{p}\right)$ is the Legendre symbol.

3) If $F$ has distinct roots in the algebraic closure $\overline{\mathbb F}_p$ and if $\overline{C}$ is the smooth projective model of  
$C$, 
then  $\overline{C}$  is an elliptic curve and satisfies
$$
|\overline{C}({\mathbb F}_p)|\equiv 1-A_{p-1}(F)\ \ \ \text{mod}\ \ \ p.
$$
\end{lemma}

{\it Proof}. Assertion 1 is proved  in \cite{silverman}, p. 141. 

Assertion 2 can be proved similarly. 

Let us prove assertion 3.
If $\deg(F)=3$ then $\overline{C}$ is a cubic in ${\mathbb P}^2$ and the complement of $C$ in $\overline{C}$ consists of one ${\mathbb F}_p$-rational point and assertion 3 follows. If $\deg(F)=4$ then, by \cite{silverman}, p. 27, $\overline{C}$ can be embedded in ${\mathbb P}^3$ such that $\overline{C}$ minus an appropriate 
 hyperplane is ${\mathbb F}_p$-isomorphic to $C$ and the intersection of $\overline{C}$ with the hyperplane consists 
of $2$ points $(0:0:\pm \sqrt{a}: 1)$. So these $2$ extra points are ${\mathbb F}_p$-rational if and only if $\left(\frac{a}{p}\right)=1$ and, again, assertion 3 follows.
\qed

\begin{corollary}
Let $A=\bZ_p$, let $c=(c_1,c_2)\in A^2$ be such that $N(c_1,c_2)\not\equiv 0$ mod $p$ (where $N$ is as in \ref{definition of N}), let $E_c$ be as in 
\ref{definition of Ec}, 
 let $F$ be as in \ref{quartic in one variable}, and let $A_{p-1}=A_{p-1}(F)$ be the Hasse invariant of $F$. Then 
 $$A_{p-1}(c_1,c_2)\equiv 0\ \ \ \text{mod}\ \ \ p$$ if and only if $E_c$ has supersingular reduction mod $p$.
\end{corollary}

{\it Proof}.
Recall from \ref{completed isogeny} that there is a separable isogeny between the smooth projective models of the reductions mod $p$ of $E_c$ and $E'_c$; so one is  supersingular if and only if the other one is. But by assertion 3 in Lemma \ref{well known} the smooth projective model of the reduction mod $p$ of $\overline{E}'_c$ is supersingular if and only if $A_{p-1}(c_1,c_2)\equiv 0$ mod $p$, hence the same is true about $E_c$. 
\qed

\begin{lemma}
\label{o lemutza}
Let $F\in A[z_1,z_2][x]$ be the quartic polynomial in \ref{quartic in one variable} and consider its Hasse invariant
$A_{p-1}:=A_{p-1}(F)\in A[z_1,z_2]$.
Then $A_{p-1}$ is homogeneous in $z_1,z_2$ of degree $\frac{p-1}{2}$ and 
\begin{equation}
\label{non congruence}
A_{p-1}\not\equiv 0\ \ \ \text{mod}\ \ \ p.\end{equation}
\end{lemma}

{\it Proof}.
Since the polynomial $F\in A[z_1,z_2][x]$ in \ref{quartic in one variable} is weighted homogeneous of degree $4$ in $z_1,z_2,x$ with respect to the weights $(2,2,1)$ we have that $F^{\frac{p-1}{2}}$ is weighed homogeneous of degree $2p-2$ with respect to the same weights. Hence $A_{p-1}(F)\in A[z_1,z_2]$ is weighted homogeneous of degree $p-1$ in $z_1,z_2$ with respect to the weights $(2,2)$, equivalently homogeneous of degree $\frac{p-1}{2}$ in $z_1,z_2$. To check \ref{non congruence}
set $t_1=z_1-a_2z_2$, $t_2=-z_1+a_1z_2$. Since 
 $$\det\left(\begin{array}{rr} 1 & -a_2\\ - 1 & a_1\end{array}\right)\not\equiv 0\ \ \ \text {mod}\ \ p,$$
 we have $A[z_1,z_2]=A[t_1,t_2]$
 and
 $$F=((a_2-a_3)x^2+t_1)((a_2-a_3)x^2+t_2).$$
  Then we have
 $$A_{p-1}=(a_2-a_3)^{\frac{p-1}{2}}t_2^{\frac{p-1}{2}}+\rho,$$
 with $\rho \in A[t_1,t_2]$ a polynomial of degree $<\frac{p-1}{2}$ in the variable $t_2$ and \ref{non congruence} follows.

\qed

\section{Statement of the main results}

In this section we introduce some more notation and state our main results.

Recall our standing assumption that $A$ a discrete valuation ring with maximal ideal generated by $p$ and residue field algebraic over ${\mathbb F}_p$. We also consider the objects 
$$H_1,H_2\in A[x_1,x_2,x_3],$$
$$F\in A[z_1,z_2][x],$$
$$N \in A[z_1,z_2],$$
$$E_c, \omega_c,$$
introduced in section \ref{section 2}. Moreover we let
$$A_{p-1}=A_{p-1}(F),\ \ \ A_{p-1}=A_{p-1}(z_1,z_2)\in A[z_1,z_2]$$
be the Hasse invariant of $F$ introduced in section 5 and we set
  \begin{equation}
  \label{definition of Q}
  Q=x_1x_2\cdot N(H_1,H_2)\cdot A_{p-1}(H_1,H_2)\in A[x_1,x_2,x_3]
  \end{equation}
   and 
   \begin{equation}
   \label{cupcake}X=Spec\ A[x_1,x_2,x_3][1/Q].\end{equation}
   In particular  $Q\not\equiv 0$ mod $p$ in $A[x_1,x_2,x_3]$ by Lemma \ref{o lemutza}.  
  
  Assume in addition that $c=(c_1,c_2)\in A^2$ satisfies
  \begin{equation}
  \label{777}
  \d c_1=\d c_2=0\ \ \ \text{and}\ \ \ N(c_1,c_2)\cdot A_{p-1}(c_1,c_2)\in A^{\times}
  \end{equation}
    and let $\d^{X}$ be any $p$-derivation on $\widehat{X}$ satisfying 
    $$\d^X H_1=\d^X H_2=0.$$
  Since $\phi^{X}(H_i)=H^p_i$ and $\phi(c_i)=c^p_i$, the Frobenius lift 
  $\phi^{X}$ on $\cO(\widehat{X})$ induces a Frobenius lift
  $\phi_c:=\phi^{E^0_c}$  on 
  $\widehat{E^0_c}$
  where $E^0_c$ is the open set of 
  $$E_c=Spec\ (A[x_1,x_2,x_3]/(H_1-c_1,H_2-c_2))$$
  given by
  $$\begin{array}{rcl}
  E^0_c & := & E_c\cap X\\
  \  &\  & \  \\
  \  & = & Spec\ (A[x_1,x_2,x_3][1/Q]/(H_1-c_1,H_2-c_2))\\
    & \  & \  \\
    \  & = & Spec\ (A[x_1,x_2,x_3][1/x_1x_2]/(H_1-c_1,H_2-c_2)).\end{array}$$
 We refer to $\phi_c$ as the {\it Frobenius lift} on $\widehat{E^0_c}$ attached to $\d^X$.
 
 Note by the way that, if $E_c^{\prime 0}$ is the open set of $E'_c$ given by
 $$E_c^{\prime 0}=Spec\ A[x,y,y^{-1}]/(y^2-F(c_1,c_2,x))$$
 then
 $$E^0_c=\pi^{-1}(E^{\prime 0}_c),$$
 where $\pi:E_c\ra E'_c$ is as in \ref{isogeny}.
   
   On the other hand,  the  global $1$-form $\omega_c$ in
    \ref{the form on intersections of two quadrics}  restricted to $E^0_c$ will be referred to as the {\it canonical} $1$-form on $E^0_c$ and will still be denoted  by $\omega_c$.
    
    The following provides  an arithmetic analogue of the classical Euler flow: 
  assertions 1 and 2   below are arithmetic analogues of assertions 1 and 2 in Theorem \ref{classical theorem} respectively.
  
  \begin{thm}
  \label{linearization theorem}
 There exists  a $p$-derivation $\d^X$ on $\widehat{X}$ such that:
 
 1) The following equalities hold in $\cO(\widehat{X})$:
  $$
   \d^XH_1=\d^XH_2=0;$$
   
   2) For
 any point $c=(c_1,c_2)\in A^2$ with
 $$
 \d c_1=\d c_2=0,\ \ \ \text{and}\ \ \ N(c_1,c_2)\cdot K(c_1,c_2)\not\equiv 0\ \ \ \text{mod}\ \ \ p,$$
  the  Frobenius lift $\phi_c$ on $\widehat{E_c^0}$ attached to $\d^X$ and the canonical $1$-form $\omega_c$ on $E_c^0$ satisfy the following congruence on $\widehat{E^0_c}$: 
  $$
  \frac{\phi_c^*}{p}\omega_c\equiv A_{p-1}(c_1,c_2)^{-1} \cdot \omega_c\ \ \ \text{mod}\ \ \ p.$$
\end{thm}

\begin{remark}
It would be interesting to see whether a strengthening of Theorem \ref{linearization theorem} holds in which the last congruence  in assertion 2 is replaced by  an equality and $A_{p-1}$ replaced by an appropriate $p$-adic modular form.
\end{remark}

\begin{definition}
A $p$-derivation $\d^X$  satisfying conditions 1 and 2 in  Theorem \ref{linearization theorem} is called an {\it arithmetic Euler flow}  (on $\widehat{X}$, attached to $H_1,H_2$). \end{definition}

  \begin{remark}\label{vin vin vin}
  The proof will show that one can find an arithmetic Euler flow $\d^X$  such that 
  $$\d^Xx_3=\frac{V}{A_{p-1}(H_1,H_2)},$$
  where $V\in A[x_1,x_2,x_3]$ is a homogeneous polynomial of degree $2p-1$. In particular, since $A_{p-1}(H_1,H_2)$ is a homogeneous polynomial of degree $p-1$ in $x_1,x_2,x_3$, we get that $\phi^X(x_3)$ is, again,  a quotient of a homogeneous polynomial of degree $2p-1$ by the homogeneous polynomial $A_{p-1}(H_1,H_2)$. Note that, for such a $\phi^X$,  the element
  $\phi^X(x_3)\in \cO(\widehat{X})$ is actually in $\cO(X)$; as we will see this automatically implies that  $\phi^X(x_1),\phi^X(x_2)$ are integral over $\cO(X)$, cf. Lemma \ref{appointment iri}. 
   \end{remark}
   
 The next result gives a parametrization of {\it arithmetic Euler flows} ``mod $p$."
 Remarkably the parameterization will involve the ``prime integrals mod $p$" of the {\it classical Euler flow}!

 Indeed let us say that two $p$-derivations  $\d^X$ and $\tilde{\d}^X$ on $\widehat{X}$ are congruent mod $p$ if 
 $$\d^X K\equiv \tilde{\d}^X K\ \ \ \text{mod}\ \ \ p$$
 for all $K\in \cO(\widehat{X})$; equivalently if the attached Frobenius lifts $\phi^X$ and $\tilde{\phi}^X$ on $\widehat{X}$ are congruent mod $p^2$ in the sense that they induce the same
 endomorphism of 
 $$\cO(\widehat{X})/p^2\cO(\widehat{X})=\cO(X)/p^2\cO(X).$$
  Also let ${\mathbb F}=A/pA$ be the residue field of $A$ and let $\overline{X}=X\otimes_A {\mathbb F}$. More generally, in what follows,  an overline will always denote a ``class mod $p$."  Next denote by
 $\overline{\d}^{cl}$  the derivation on $\cO(\overline{X})$ obtained by reducing mod $p$ the classical Euler flow, i.e.,
 \medskip
 \begin{equation}
 \begin{array}{rcl}
 \overline{\d}^{cl} x_1 & = & (a_2-a_3)x_2x_3,\\
  \overline{\d}^{cl} x_2 & = & (a_3-a_1)x_3x_1,\\
   \overline{\d}^{cl} x_3 & = & (a_1-a_2)x_1x_2.\end{array}\end{equation}
   \medskip
   The elements $\overline{K}\in \cO(\overline{X})$ such that 
  $$\overline{\d}^{cl}\overline{K}=0$$
  can be referred to as the {\it prime integrals mod $p$} of the classical Euler flow. 
  
  \medskip
  
Then we have:

\begin{thm}
\label{principal}
Assume ${\mathbb F}$ is infinite. Then the quotient set
$$\frac{\{\text{arithmetic Euler flows on $\widehat{X}$}\}}{\{\text{congruence mod $p$}\}}$$
is a  principal homogeneous space for the ${\mathbb F}$-linear space
$$\{\overline{K}\in \cO(\overline{X}); \  \overline{\d}^{cl}\overline{K}=0\}.$$
\end{thm}

\begin{remark}
Recall that a set is a principal homogeneous space for a group if it has a simply transitive action of that group. A principal homogeneous space is either empty or in bijection with the group. The principal homogeneous space in Theorem \ref{principal} is non-empty due to Theorem \ref{linearization theorem}.
\end{remark}

\begin{remark}
The ${\mathbb F}$-linear space in Theorem \ref{principal} is, of course, rather large:
it contains the ${\mathbb F}$-subalgebra of $\cO(\overline{X})$ generated by
$$\overline{H}_1,\ \ \ \overline{H}_2,\ \ \ \overline{A}_{p-1}(\overline{H}_1,\overline{H}_2)^{-1},\ \ \ \overline{N}(\overline{H}_1,\overline{H}_2)^{-1},$$
 and by all the $p$-powers of elements of $\cO(\overline{X})$.
\end{remark}
 
 \section{Proof of the main results}

  We first concentrate on proving  our Theorem \ref{linearization theorem}; we need a preparation. 
  
  For the discussion below $Q$ can be any element
  of $A[x_1,x_2,x_3]$ with $Q\not\equiv 0$ mod $p$ (which is not necessarily as in \ref{definition of Q}) and we continue to denote by $X$ the scheme in \ref{cupcake}.
  Consider the morphism
$${\mathcal H}:X\ra {\mathbb A}^2,$$
obtained from \ref{condition 4} by restriction, hence given by
$ z_i\mapsto H_i$, $i=1,2,$
and continue to denote by 
$${\mathcal H}:\widehat{X}\ra \widehat{{\mathbb A}^2}$$
the induced morphism. 
We will also be considering the {\it trivial}  Frobenius lift $\phi_0^{{\mathbb A}^2}$ on $\widehat{{\mathbb A}^2}$
given by 
$\phi_0^{{\mathbb A}^2}(z_i)=z_i^p$. Its attached $p$-derivation satisfies $\d_0^{{\mathbb A}^2} z_i=0$. 

\begin{lemma}
\label{appointment iri} 
Assume $Q$ is in the ideal  $x_1x_2A[x_1,x_2,x_3]$.
Let ${\mathcal F}$ be the set of  Frobenius lifts $\phi^{X}$ on $\widehat{X}$ such that the following diagram is commutative:
\begin{equation}
\label{divine comedy bleep}
\begin{array}{rcl}
\widehat{X} & \stackrel{\phi^{X}}{\longrightarrow} & \widehat{X}\\
{\mathcal H} \downarrow & \  & \downarrow {\mathcal H}\\
\widehat{{\mathbb A}^2} & \stackrel{\phi_0^{{\mathbb A}^2}}{\longrightarrow} & \widehat{{\mathbb A}^2}
\end{array}
\end{equation}
Then the following hold:

1) The set ${\mathcal F}$  is in  bijection with  the ring $\cO(\widehat{X})$; the bijection is given by 
\begin{equation}
\label{bijection breast}
{\mathcal F}\ra \cO(\widehat{X}),\ \ \ \phi^X\mapsto \Delta_3,\ \ \ \text{where}\ \ \ \phi^X(x_3)=x^p_3+p\Delta_3.
\end{equation}

2) If $\phi^X_3(x_3)$ (equivalently $\Delta_3$)  is in $\cO(X)$  then $\phi^X(x_1)$ and $\phi^X(x_2)$ are integral over $\cO(X)$.

3) Two Frobenius lifts $\phi^X$ and $\tilde{\phi}^X$ in ${\mathcal F}$
are congruent mod $p^2$ if and only if the corresponding elements $\Delta_3$ and $\tilde{\Delta}_3$ are congruent mod $p$ in the ring $\cO(\widehat{X})$.
\end{lemma}

\begin{remark}
Let $\d^{X}$ be the $p$-derivation attached to $\phi^{X}$. Then  the commutativity of diagram \ref{divine comedy bleep} is equivalent to the condition:
$$\d^{X} H_i=0,\ \ \ i=1,2.$$
\end{remark}

\medskip

\begin{remark}
We will give two proofs of assertion 1 in Lemma \ref{appointment iri}. The first one is more conceptual  but non-explicit. The second one is more computational but also more explicit, and 
leads to proofs of assertions 2 and 3 as well. Assertion 2 can be used to prove the existence of {\it correspondence structures} 
(in the sense of \cite{curvature2}) for our Frobenius lifts $\phi^X$ in case $\Delta_3$ is  in $\cO(X)$. However we will not be concerned here with these correspondence structures.
\end{remark}

For the first proof of assertion 1 in Lemma \ref{appointment iri} 
we need to recall the following:

\begin{lemma} 
\label{etale etale etale}
Let $X\ra Y$ be an \'{e}tale morphism of affine schemes of  finite type over $A$. Then any $p$-derivation  $\cO(\widehat{Y})\ra \cO(\widehat{X})$ lifts to a unique $p$-derivation $\cO(\widehat{X})\ra \cO(\widehat{X})$.
\end{lemma}

{\it Proof}. Cf. \cite{book}, Lemma 3.14, p. 76.\qed

\medskip

{\it First proof of assertion 1 in Lemma \ref{appointment iri}}.
Let $Y={\mathbb A}^3=Spec\ A[z_1,z_2,z_3]$ be an affine $3$-space with coordinates $z_1,z_2,z_3$.
and consider the morphism 
$$\tilde{\mathcal H}:X\ra Y$$
 given by
$$z_1\mapsto H_1,\ \ z_2\mapsto H_2,\ \ z_3\mapsto x_3.$$
Clearly $\tilde{\mathcal H}$ composed with the projection $Y={\mathbb A}^3\ra {\mathbb A}^2$ onto the first $2$ components equals ${\mathcal H}$.
Let $\Delta_3\in \c(\widehat{X})$ be an arbitrary element; we shall construct a Frobenius lift  $\phi^X\in {\mathcal F}$ such that $\phi^X(x_3)=x_3^p+p\Delta_3$ and this will end the proof of assertion 1. Note that the Jacobian matrix of $\tilde{\mathcal H}$ equals
$$\left(\begin{array}{ccc} 2a_1x_1 & 2a_2x_2 & 2a_3x_3\\
\  & \  & \  \\
2x_1 & 2x_2 & 2x_3\\
\  & \  & \  \\
0& 0 & 1\end{array}\right)$$
so $\tilde{\mathcal H}$ is \'{e}tale because $x_1x_2$ is invertible on $X$.
Consider the $p$-derivation $\d^{YX}:\cO(\widehat{Y})\ra \cO(\widehat{X})$
 defined by
$$\d^{YX}(z_1)=0,\ \ \ \d^{YX}(z_2)=0,\ \ \ \d^{YX}(z_3)=\Delta_3.$$
By Lemma \ref{etale etale etale} $\d^{YX}$ extends to a unique $p$-derivation
$\d^X:\cO(\widehat{X})\ra \cO(\widehat{X})$, hence we will have
$$\d^X(H_1)=0,\ \ \ \d^X(H_2)=0,\ \ \ \d^X(x_3)=\Delta_3,$$
so the Frobenius lift attached to $\d^X$ is in ${\mathcal F}$ and satisfies the desired requirement.
\qed

\medskip

{\it Second proof of assertion 1 Lemma \ref{appointment iri}}. 
Let $\phi^{X}$ be a Frobenius lift on $\widehat{X}$ and let
$\Delta_i\in \cO(\widehat{X})$  
be such that
$$\Phi_i:=x_i^p+p\Delta_i
=\phi^{X}(x_i),\ \ \ i=1,2,3.$$
For the sake of uniformity in notation write $a_{1j}=a_j$ and $a_{2j}=1$.
Then the commutativity of the diagram 
\ref{divine comedy bleep}
 is equivalent to the condition
 \begin{equation}
\label{cape cod}
\sum_{j=1}^3\phi(a_{ij})\Phi_j^2=(\sum_{j=1}^3a_{ij}x_j^2)^p,\ \ \ i=1,2.
\end{equation}
We can rewrite the system \ref{cape cod} as
 \begin{equation}
\label{cape cod 1}
\sum_{j=1}^2\phi(a_{ij})\Phi_j^2=(\sum_{j=1}^3a_{ij}x_j^2)^p-
\phi(a_{i3})(x_3^p+p\Delta_3)^2,\ \ \ i=1,2.
\end{equation}
Note that the right hand side of \ref{cape cod 1} is congruent mod $p$ to
$$\sum_{i=1}^2 a_{ij}^px_j^{2p}.
$$
Fix an element $\Delta_3$ in $\cO(\widehat{X})$ and let $M$ be the matrix $(a_{ij})_{i,j=1}^2$.
By Cramer's rule the  system \ref{cape cod 1}, viewed as a linear system with unknowns  $\Phi_1^2,\Phi_2^2$, has a unique solution given by
$$
\Phi_i^2=\frac{D_i}{D}
$$
where $D=\det(\phi(M))$ and $D_i$ is the determinant of the matrix obtained from the matrix $\phi(M)$ by replacing the $i$-th column with the column vector appearing in the right hand side of the system \ref{cape cod 1}.
Now clearly
\medskip
$$D \equiv \det\left(\begin{array}{cc}a_{11}^p &  a_{12}^p\\
\  & \  \\
a_{21}^p &  a_{2m}^p\end{array}\right)\ \ \ \text{mod}\ \ p.$$
\medskip
On the other hand,
\medskip
$$D_1 \equiv \det\left(\begin{array}{cc} \sum_{j=1}^2 a_{1j}^px_j^{2p}  & a_{12}^p\\
\  & \   \\
 \sum_{j=1}^2 a_{2j}^px_j^{2p} &  a_{22}^p\end{array}\right)=x_1^{2p}\cdot D\ \ \ \text{mod}\ \ p,$$
 \medskip
 and, similarly,
 \medskip
 $$D_2\equiv x_2^{2p}\cdot D\ \ \ \text{mod}\ \ p.$$
\medskip
Hence we have
$$\Phi_i^2=x_i^{2p}(1+pG_i),\ \ \ i=1,2,$$
for some $G_i\in \cO(\widehat{X})$. Hence the system \ref{cape cod 1}, viewed as a non-linear system with unknowns $\Phi_1,\Phi_2$, has (under the restriction $\Phi_i\equiv x_i^p$ mod $p$) a unique solution given by
\begin{equation}
\Phi_i=x_i^p(1+pG_i)^{1/2}\in \cO(\widehat{X}),
\end{equation}
where the $1/2$ power is the principal root (the root congruent to $1$ mod $p$). This ends our proof of assertion 1. \qed

\medskip

{\it Proof of assertions 2 and 3 in Lemma \ref{appointment iri}}.
To check assertion 2 note that if $\Delta_3$  in the second proof above is in $\cO(X)$ then so are $D_1,D_2$, and hence so are $\Phi^2_1,\Phi^2_2$; so $\phi^X(x_1),\phi^X(x_2)$ are integral over $\cO(X)$.

The ``only if" part of assertion 3 is clear. To check the ``if" part assume $\Delta_3\equiv \tilde{\Delta}_3$ mod $p$ for two Frobenius lifts $\phi^X,\tilde{\phi}^X\in {\mathcal F}$. Then, clearly $\phi^X(x_3)\equiv \tilde{\phi}^X(x_3)$ mod $p^2$ and, by the second proof above, 
$$p^2\ |\ \Phi_1^2- \tilde{\Phi}_1^2=(\Phi_1-\tilde{\Phi}_1)(\Phi_1+\tilde{\Phi}_1),$$
where $\tilde{\Phi}_i=\tilde{\phi}^X(x_i)$.
Since $p$ is prime in $\cO(\widehat{X})$ and $p$ does not divide $\Phi_1+\tilde{\Phi}_1$ (the latter being congruent to $2x_1^p$ mod $p$) 
it follows that 
$$p^2\ |\ \Phi_1-\tilde{\Phi}_1$$
in $\cO(\widehat{X})$ so $\phi^X(x_1)\equiv \tilde{\phi}^X(x_1)$ mod $p^2$. Similarly $\phi^X(x_2)\equiv \tilde{\phi}^X(x_2)$ mod $p^2$ and we are done.
\qed

  \bigskip
  
  At this point we are ready to conclude the proof of our Theorem \ref{linearization theorem}.
  
  \bigskip
 
{\it Proof of Theorem \ref{linearization theorem}}. 
Write 
 \begin{equation}
 \label{one star}
 A_{p-1}(z_1,z_2)^{-1}\cdot F(x,z_1,z_2)^{\frac{p-1}{2}}-x^{p-1}=\sum_{i=0}^{2p-2} R_i(z_1,z_2)x^i,
 \end{equation}
 with
 $$R_i(z_1,z_2)\in A_{p-1}(z_1,z_2)^{-1}\cdot A[z_1,z_2].$$
 Since $R_{p-1}=0$ we have that
 $$S(z_1,z_2,x):=\sum_{i=0}^{2p-2} R_i(z_1,z_2)\cdot \frac{x^{i+1}}{i+1}\in 
 A_{p-1}(z_1,z_2)^{-1}\cdot A[z_1,z_2][x].$$
 Define 
 $$\Delta_3:=S(H_1,H_2,x_3)\in  A_{p-1}(H_1,H_2)^{-1}\cdot A[x_1,x_2,x_3]\subset \cO(X).$$
 By the bijection \ref{bijection breast} in Lemma  \ref{appointment iri}
 there is a $p$-derivation $\d^X$ on $\widehat{X}$ killing $H_1$ and $H_2$ whose Frobenius lift $\phi^X$ satisfies
 $$\phi^X(x_3)=x_3^p+p\Delta_3.$$
 Let now $(c_1,c_2)\in A^2$ satisfy the conditions in assertion 2 of Theorem \ref{linearization theorem}.  We claim that the congruence in the conclusion of assertion 2 holds and this will end the proof.
 Indeed we have 
 \medskip
 $$
 \frac{\phi_c^*}{p}\omega_c = \frac{\phi_c^*}{p}\left( i_c^* \frac{dx_3}{(a_1-a_2)x_1x_2}\right)=  \frac{(d\Phi_{3c})/p}{\phi(a_1-a_2) \Phi_{1c}\Phi_{2c}},
 $$
 \medskip
 where, for $i=1,2,3$,   $\Phi_{ic}=i_c^* \Phi_i$ is the restriction of $\Phi_i=\phi^X(x_i)$ to $E^0_c$.
 Now
 $$
 \begin{array}{rcl}
 \Phi_{3c} & = & x_3^p+pi^*_c\Delta_3\\
 \  & \  & \ \\
 \  & = & x_3^p+ pS(c_1,c_2,x_3)\\
 \  & \  & \  \\
 \  & = & x_3^p+p\sum_{i=0}^{2p-2} R_i(c_1,c_2) \cdot \frac{x_3^{i+1}}{i+1},
 \end{array}
 $$
 where, in order to simplify the notation, we continued to write $x_i$ in place of $i^*_cx_i$. 
 Hence, with the same abuse of notation,
 $$
 \begin{array}{rcll}
 (d\Phi_{3c})/p & = & (x_3^{p-1}+\sum_{i=0}^{2p-2}R_i(c_1,c_2)x_3^i)dx_3 & \  \\
 \  & \  & \  & \  \\
 \  & = & A_{p-1}(c_1,c_2)^{-1}F(x_3,c_1,c_2)^{\frac{p-1}{2}}\cdot dx_3 & (\text{by}\ \ \ref{one star})\\
 \  & \  & \  & \  \\
 \  & = & A_{p-1}(c_1,c_2)^{-1}((a_2-a_3)x_3^2+c_1-a_2c_2)^{\frac{p-1}{2}} \times & \ \\
  \  & \  & \  & \  \\
 \  & \  & \times 
 ((a_3-a_1)x_3^2-c_1+a_1c_2)^{\frac{p-1}{2}}dx_3 &\ \\
   \  & \  & \  & \  \\
   \  & = & A_{p-1}(c_1,c_2)^{-1} (a_1-a_2)^{p-1}x_1^{p-1}x_2^{p-1} dx_3 & (\text{by}\ \  \ref{congruences mod I}).
 \end{array}
 $$
 So we have
 $$
 \begin{array}{rcll}
 \frac{\phi_c^*}{p}\omega_c & \equiv & A_{p-1}(c_1,c_2)^{-1} \cdot \frac{(a_1-a_2)^{p-1}x_1^{p-1}x_2^{p-1}dx_3}{(a_1-a_2)^px_1^px_2^p} & \text{mod}\ \ \ p\\
 \  & \  & \  & \  \\
 \  & \equiv & A_{p-1}(c_1,c_2)^{-1} \frac{dx_3}{(a_1-a_2)x_1x_2} & \text{mod}\ \ \ p\\
 \  & \  & \  & \  \\
\  & = & A_{p-1}(c_1,c_2)^{-1}\cdot \omega_c.
 \end{array}
 $$
 This ends our proof.
 \qed
 
 \begin{remark}
 We claim that the arithmetic Euler flow $\d^X$ constructed in the proof above satisfies
 $$\d^X x_3=\frac{V}{A_{p-1}(H_1,H_2)}$$
 with $V\in A[x_1,x_2,x_3]$ homogeneous in $x_1,x_2,x_3$ of degree $2p-1$.
 To check this recall that
  the Hasse invariant $A_{p-1}=A_{p-1}(F)\in A[z_1,z_2]$ is weighted homogeneous of degree $p-1$ in $z_1,z_2$ with respect to the weights $(2,2)$ so $A_{p-1}(H_1,H_2)\in A[x_1,x_2,x_3]$ is homogeneous in $x_1,x_2,x_3$ of degree $p-1$. On the other hand the elements 
 $$A_{p-1}(z_1,z_2) \cdot R_i(z_1,z_2)$$ in our proof are 
 weighted homogeneous polynomials of degree $2p-2-i$ in $z_1,z_2$ with respect to the weights $(2,2)$,
 hence
 $$A_{p-1}\cdot S(z_1,z_2,x)$$
 is a weighted homogeneous polynomial of degree $2p-1$ in $z_1,z_2,x$ with respect to the weights $(2,2,1)$, 
  and hence
 $$V(x_1,x_2,x_3):=A_{p-1}(H_1,H_2)\cdot \Delta_3,$$
  is a homogeneous polynomial in $A[x_1,x_2,x_3]$ of degree $2p-1$ in $x_1,x_2,x_3$. \end{remark}
 
  \medskip
  
  We next concentrate on proving Theorem \ref{principal}.
  
  We need the following Lemma:
  
  \begin{lemma}
  \label{singe}
 Let $\overline{K}\in \cO(\overline{X})$ and $c=(c_1,c_2)\in A^2$ with 
 $\d c_1=\d c_2=0$, $N(c_1,c_2)\not\equiv 0$ mod $p$. Let $i_c:E^0_c\ra X$ be, as usual,  the inclusion, let $\overline{i}_c:\overline{E}^0_c\ra \overline{X}$ be its reduction mod $p$, and consider the function $\overline{i}^*_c\overline{K}\in \cO(\overline{E}^0_c)$. 
 Let $\overline{\omega}_c$ be the $1$-form on $\overline{E}^0_c$ obtained from $\omega_c$ by reduction mod $p$. 
 Then we have the following equality of $1$-forms on $\overline{E}^0_c$:
 $$d(\overline{i}_c^*\overline{K})=\overline{i}_c^*(\overline{\d}^{cl}\overline{K})\cdot \overline{\omega}_c.$$
  \end{lemma}
  
  {\it Proof}.
  We have the following computation (in which ``$\sum_{(123)}$" means ``sum obtained by cyclic permutations of the indices $1,2,3$"):
  $$
  \begin{array}{rcl}
  d(\overline{i}_c^*\overline{K}) & = & \overline{i}_c^*(d\overline{K})\\
  \  & \  & \  \\
  \  & = & \overline{i}_c^*(\frac{\partial \overline{K}}{\partial x_1}dx_1+
  \frac{\partial \overline{K}}{\partial x_2}dx_2+\frac{\partial \overline{K}}{\partial x_3}dx_3)\\
  \  & \  & \  \\
  \  & = & \overline{i}_c^*(\sum_{(123)} (a_2-a_3)x_2x_3\frac{\partial \overline{K}}{\partial x_1}\frac{dx_1}{(a_2-a_3)x_2x_3})\\
   \  & \  & \  \\
  \  & = & \sum_{(123)} \overline{i}_c^*((a_2-a_3)x_2x_3\frac{\partial \overline{K}}{\partial x_1})\cdot \overline{i}_c^*\frac{dx_1}{(a_2-a_3)x_2x_3}\\
   \  & \  & \  \\
  \  & = & (\sum_{(123)} \overline{i}_c^*((a_2-a_3)x_2x_3\frac{\partial \overline{K}}{\partial x_1}))\cdot \overline{\omega}_c\\
    \  & \  & \  \\
  \  & = &
\overline{i}_c^*(\overline{\d}^{cl}\overline{K})\cdot \overline{\omega}_c.
    \end{array}
  $$
  \qed
  
  \medskip

  {\it Proof of Theorem \ref{principal}}.
  By Lemma \ref{appointment iri} the set of arithmetic Euler flows $\d^X$  modulo $p$
  is in a natural bijection with the set of all $\overline{\Delta}_3$ in $\cO(\overline{X})$ satisfying 
 the following equality of $1$-forms 
 \begin{equation}
 \label{vomit}
 \frac{x_3^{p-1}dx_3+d(\overline{i}^*_c\overline{\Delta}_3)}{(a_1-a_2)^px_1^px_2^p}=A_{p-1}(c_1,c_2)^{-1}\frac{dx_3}{(a_1-a_2)x_2x_3}
 \end{equation}
 on $\overline{E}^0_c$ for all $c=(c_1,c_2)\in A^2$ with 
 \begin{equation}
 \label{mucus}
 \d c_1=\d c_2=0,\ \ \ N(c_1,c_2)A_{p-1}(c_1,c_2)\not\equiv 0\ \ \ \text{mod}\ \ \ p.\end{equation}
  The bijection is given by the rule $\phi^X(x_3)=x_3^p+p\Delta_3$.
  So if $\tilde{\phi}^X=x_3^p+p\tilde{\Delta}_3$ is another arithmetic Euler flow and we set $K=\Delta_3-\tilde{\Delta}_3\in \cO(\widehat{X})$ then its redulction mod $p$, $\overline{K}\in \cO(\overline{X})$, satisfies
  $$d(\overline{i}_c^*\overline{K})=0$$
  for all $c$ that satisfy \ref{mucus}. By Lemma \ref{singe} we have that $\overline{K}$ satisfies
  $$\overline{i}_c^*(\overline{\d}^{cl}\overline{K})=0$$
  for all $c$ that satisfy \ref{mucus}. Now, for any point in ${\mathbb F}^2$ which is not a zero of $N\cdot A_{p-1}$, the Teichm\"{u}ller lift $c$ of that point satisfies \ref{mucus}. Since ${\mathbb F}$ is infinite the set of points in ${\mathbb F}^2$ which are not zeros of  $N\cdot A_{p-1}$ is Zariski dense in ${\mathbb A}^2_{\mathbb F}$ hence the
   union of all curves $\overline{E}^0_c$ with $c$ as in \ref{mucus} is Zariski dense in 
  $\overline{X}$. We conclude that $\overline{\d}^{cl}\overline{K}=0$.
  
  Conversely, if $\overline{\Delta}_3$ satisfies \ref{vomit} for all $c$ satisfying \ref{mucus} and if $\overline{\d}^{cl}\overline{K}=0$ then $\overline{\Delta}_3+\overline{K}$ will still satisfy \ref{vomit} for all $c$ satisfying \ref{mucus}. This shows that the ${\mathbb F}$-linear space in the statement of our theorem acts on the set of arithmetic Euler flows  mod $p$; by the first part of the argument this action is  transitive. Also  the isotropy groups are clearly trivial. This ends our proof.
    \qed

\end{document}